\documentclass[11pt]{amsart}
\usepackage{amssymb}
\usepackage{amsmath}
\usepackage[active]{srcltx}
\usepackage{t1enc}
\usepackage[latin2]{inputenc}
\usepackage{verbatim}
\usepackage{amsmath,amsfonts,amssymb,amsthm}
\usepackage[mathcal]{eucal}
\usepackage{enumerate}
\usepackage[centertags]{amsmath}
\usepackage{graphics}

\newtheorem{theorem}{Theorem}

\begin{document}
\author{George Tephnadze}
\title[Fourier coefficients and partial sums]{A note on the Fourier
coefficients and partial sums of Vilenkin-Fourier series}
\address{G. Tephnadze, Department of Mathematics, Faculty of Exact and
Natural Sciences, Tbilisi State University, Chavchavadze str. 1, Tbilisi
0128, Georgia}
\email{giorgitephnadze@gmail.com}
\date{}
\maketitle

\begin{abstract}
The main aim of this paper is to investigate Paley type and Hardy-Littlewood
type inequalities and strong convergence theorem of partial sums of
Vilenkin-Fourier series.
\end{abstract}

\date{}



Let $\mathbb{N}_{+}$ denote the set of the positive integers, $\mathbb{N}:=%
\mathbb{N}_{+}\cup \{0\}.$

Let $m:=(m_{0},m_{1},...)$ denote a sequence of the positive numbers, not
less than 2.

Denote by
\begin{equation*}
Z_{m_{k}}:=\{0,1,...,m_{k}-1\}
\end{equation*}
the additive group of integers modulo $m_{k}.$

Define the group $G_{m}$ as the complete direct product of the group $%
Z_{m_{j}}$ with the product of the discrete topologies of $Z_{m_{j}}`$s.

The direct product $\mu ,$ of the measures
\begin{equation*}
\mu _{k}\left( \{j\}\right) :=1/m_{k},\text{\qquad }(j\in Z_{m_{k}})
\end{equation*}
is the Haar measure on $G_{m},$ with $\mu \left( G_{m}\right) =1.$

If $\sup_{n}m_{n}<\infty $, then we call $G_{m}$ a bounded Vilenkin group.
If the generating sequence $m$ is not bounded then $G_{m}$ is said to be an
unbounded Vilenkin group. \textbf{In this paper we discuss bounded Vilenkin
groups only.}

The elements of $G_{m}$ represented by sequences
\begin{equation*}
x:=(x_{0},x_{1},...,x_{j},...),\qquad \left( \text{ }x_{k}\in
Z_{m_{k}}\right) .
\end{equation*}

It is easy to give a base for the neighborhood of $G_{m}:$
\begin{equation*}
I_{0}\left( x\right) :=G_{m},
\end{equation*}%
\begin{equation*}
I_{n}(x):=\{y\in G_{m}\mid y_{0}=x_{0},...,y_{n-1}=x_{n-1}\},\text{ }(x\in
G_{m},\text{ }n\in \mathbb{N}).
\end{equation*}%
Denote $I_{n}:=I_{n}\left( 0\right) ,$ for $n\in \mathbb{N}$ and $\overset{-}%
{I_{n}}:=G_{m}$ $\backslash $ $I_{n}$ .

If we define the so-called generalized number system, based on $m$ in the
following way :
\begin{equation*}
M_{0}:=1,\text{ \qquad }M_{k+1}:=m_{k}M_{k\text{ }}\ \qquad (k\in \mathbb{N}%
),
\end{equation*}%
then every $n\in \mathbb{N}$ can be uniquely expressed as $n=\overset{\infty
}{\underset{j=0}{\sum }}n_{j}M_{j},$ where $n_{j}\in Z_{m_{j}}$ $~(j\in
\mathbb{N})$ and only a finite number of $n_{j}`$s differ from zero.

Let $\left\vert n\right\vert :=\max $ $\{j\in \mathbb{N}:$ $n_{j}\neq 0\}.$

Denote by $\mathbb{N}_{n_{0}}$ the subset of positive integers $\mathbb{N}%
_{+},$ for which $n_{\left\vert n\right\vert }=n_{0}=1.$ Then for every $%
n\in \mathbb{N}_{n_{0}},$ $M_{k}<n<$ $M_{k+1}$ can be written as $%
n=M_{0}+\sum_{j=1}^{k-1}n_{j}M_{j}+M_{k}=1+\sum_{j=1}^{k-1}n_{j}M_{j}+M_{k},$
where $n_{j}\in \left\{ 0,m_{j}-1\right\} ,$ $~(j\in \mathbb{N}_{+}).$

By simple calculation we get
\begin{equation}
\underset{\left\{ n:M_{k}\leq n\leq M_{k+1},\text{ }n\in \mathbb{N}%
_{n_{0}}\right\} }{\sum }1=\frac{M_{k-1}}{m_{0}}\geq cM_{k},  \label{1a}
\end{equation}%
where $c$ is absolute constant.

Denote by $L_{1}\left( G_{m}\right) $ the usual (one dimensional) Lebesgue
space.

Next, we introduce on $G_{m}$ an ortonormal system, which is called the
Vilenkin system.

At first define the complex valued function $r_{k}\left( x\right)
:G_{m}\rightarrow
%TCIMACRO{\U{2102} }%
%BeginExpansion
\mathbb{C}
%EndExpansion
,$ the generalized Rademacher functions as
\begin{equation*}
r_{k}\left( x\right) :=\exp \left( 2\pi \iota x_{k}/m_{k}\right) ,\text{
\qquad }\left( \iota ^{2}=-1,\text{ }x\in G_{m},\text{ }k\in \mathbb{N}%
\right) .
\end{equation*}

Now define the Vilenkin system $\psi :=(\psi _{n}:n\in \mathbb{N})$ on $%
G_{m} $ as:
\begin{equation*}
\psi _{n}(x):=\overset{\infty }{\underset{k=0}{\Pi }}r_{k}^{n_{k}}\left(
x\right) ,\text{ \qquad }\left( n\in \mathbb{N}\right) .
\end{equation*}

Specifically, we call this system the Walsh-Paley one if $m\equiv 2$.

The Vilenkin system is ortonormal and complete in $L_{2}\left( G_{m}\right) $
$\,$\cite{AVD,Vi}.

Now we introduce analogues of the usual definitions in Fourier-analysis. If $%
f\in L_{1}\left( G_{m}\right) $ we can establish the the Fourier
coefficients, the partial sums, the Dirichlet kernels, with respect to the
Vilenkin system in the usual manner:
\begin{eqnarray*}
\widehat{f}\left( k\right) &:&=\int_{G_{m}}f\overline{\psi }_{k}d\mu \text{%
\thinspace \qquad\ \ }\left( \text{ }k\in \mathbb{N}\text{ }\right) , \\
S_{n}f &:&=\sum_{k=0}^{n-1}\widehat{f}\left( k\right) \psi _{k}\text{ \qquad
}\left( \text{ }n\in \mathbb{N}_{+},\text{ }S_{0}f:=0\right) , \\
D_{n} &:&=\sum_{k=0}^{n-1}\psi _{k}\text{ \qquad\ \ \ \ \ \ \ \ }\left(
\text{ }n\in \mathbb{N}_{+}\text{ }\right) ,
\end{eqnarray*}

Recall that
\begin{equation}
\quad \hspace*{0in}D_{M_{n}}\left( x\right) =\left\{
\begin{array}{l}
\text{ }M_{n}\text{ \ \ \thinspace \thinspace \thinspace \thinspace
if\thinspace \thinspace }x\in I_{n} \\
\text{ }0\ \ \ \ \ \text{\thinspace \thinspace \thinspace \thinspace
\thinspace if \thinspace \thinspace }x\notin I_{n}%
\end{array}%
\right.  \label{3}
\end{equation}%
\vspace{0pt}and
\begin{equation}
D_{n}=\psi _{n}\left( \sum_{j=0}^{\infty
}D_{M_{j}}\sum_{u=m_{j}-n_{j}}^{m_{j}-1}r_{j}^{u}\right) .  \label{3aa}
\end{equation}%
The norm (or quasinorm) of the space $L_{p}(G_{m})$ is defined by \qquad

\begin{equation*}
\left\Vert f\right\Vert _{p}:=\left( \int_{G_{m}}\left\vert f\right\vert
^{p}d\mu \right) ^{1/p}\qquad \left( 0<p<\infty \right) .
\end{equation*}%
The space $L_{p,\infty }\left( G_{m}\right) $ consists of all measurable
functions $f$ for which

\begin{equation*}
\left\Vert f\right\Vert _{L_{p,\infty }}:=\underset{\lambda >0}{\sup }%
\lambda \mu \left( f>\lambda \right) ^{1/p}<+\infty .
\end{equation*}

The $\sigma $-algebra, generated by the intervals $\left\{ I_{n}\left(
x\right) :x\in G_{m}\right\} $ will be denoted by $\digamma _{n}$ $\left(
n\in \mathbb{N}\right) .$ The conditional expectation operators relative to $%
\digamma _{n}\left( n\in \mathbb{N}\right) $ are denoted by $E_{n}.$ Then

\begin{eqnarray*}
E_{n}f\left( x\right) &=&S_{M_{n}}f\left( x\right) =\sum_{k=0}^{M_{n}-1}%
\widehat{f}\left( k\right) w_{k} \\
&=&\frac{1}{\left| I_{n}\left( x\right) \right| }\int_{I_{n}\left( x\right)
}f(x)d\mu (x),
\end{eqnarray*}
where $\left| I_{n}\left( x\right) \right| =M_{n}^{-1}$ denotes the length
of $I_{n}\left( x\right) .$

A sequence $f=\left( f^{\left( n\right) },\text{ }n\in \mathbb{N}\right) $
of functions $f_{n}\in L_{1}\left( G\right) $ is said to be a dyadic
martingale if (for details see e.g. \cite{We1})

$\left( i\right) $ $f^{\left( n\right) }$ is $\digamma _{n}$ measurable, for
all $n\in \mathbb{N},$

$\left( ii\right) $ $E_{n}f^{\left( m\right) }=f^{\left( n\right) },$ for
all $n\leq m.$

The maximal function of a martingale $f$ is denoted by \qquad
\begin{equation*}
f^{\ast }=\sup_{n\in \mathbb{N}}\left\vert f^{\left( n\right) }\right\vert .
\end{equation*}

In case $f\in L_{1},$ the maximal functions are also be given by
\begin{equation*}
f^{\ast }\left( x\right) =\sup_{n\in \mathbb{N}}\frac{1}{\left\vert
I_{n}\left( x\right) \right\vert }\left\vert \int_{I_{n}\left( x\right)
}f\left( u\right) \mu \left( u\right) \right\vert .
\end{equation*}

For $0<p<\infty ,$ the Hardy martingale spaces $H_{p}$ $\left( G_{m}\right) $
consist of all martingales, for which
\begin{equation*}
\left\| f\right\| _{H_{p}}:=\left\| f^{*}\right\| _{L_{p}}<\infty .
\end{equation*}

If $f\in L_{1},$ then it is easy to show that the sequence $\left(
S_{M_{n}}f:n\in \mathbb{N}\right) $ is a martingale. If $f=\left( f^{\left(
n\right) },n\in \mathbb{N}\right) $ is martingale, then the Vilenkin-Fourier
coefficients must be defined in a slightly different manner: $\qquad \qquad $
\begin{equation}
\widehat{f}\left( i\right) :=\lim_{k\rightarrow \infty
}\int_{G_{m}}f^{\left( k\right) }\left( x\right) \overline{\Psi }_{i}\left(
x\right) d\mu \left( x\right) .  \label{3a}
\end{equation}%
\qquad \qquad \qquad \qquad

The Vilenkin-Fourier coefficients of $f\in L_{1}\left( G_{m}\right) $ are
the same as those of the martingale $\left( S_{M_{n}}f:n\in \mathbb{N}%
\right) $ obtained from $f$ .

A bounded measurable function $a$ is p-atom, if there exist a dyadic
interval $I$, such that \qquad
\begin{equation*}
\left\{
\begin{array}{l}
a)\qquad \int_{I}ad\mu =0 \\
b)\ \qquad \left\Vert a\right\Vert _{\infty }\leq \mu \left( I\right) ^{-1/p}
\\
c)\qquad \ \text{supp}\left( a\right) \subset I.\qquad%
\end{array}%
\right.
\end{equation*}

The Hardy martingale spaces $H_{p}$ $\left( G_{m}\right) ,$ for $0<p\leq 1$
have an atomic characterization. Namely, the following theorem is true (see
\cite{We3}):

\textbf{Theorem W.} A martingale $f=\left( f^{\left( n\right) },n\in \mathbb{%
N}\right) $ is in $H_{p}\left( 0<p\leq 1\right) $ if and only if there exist
a sequence $\left( a_{k},k\in \mathbb{N}\right) $ of p-atoms and a sequence $%
\left( \mu _{k},k\in \mathbb{N}\right) $ of a real numbers, such that for
every $n\in \mathbb{N}:$

\begin{equation}
\qquad \sum_{k=0}^{\infty }\mu _{k}S_{M_{n}}a_{k}=f^{\left( n\right) },
\label{1Aaaa}
\end{equation}

\begin{equation*}
\qquad \sum_{k=0}^{\infty }\left| \mu _{k}\right| ^{p}<\infty .
\end{equation*}
Moreover, $\left\| f\right\| _{H_{p}}\backsim \inf \left( \sum_{k=0}^{\infty
}\left| \mu _{k}\right| ^{p}\right) ^{1/p}$, where the infimum is taken over
all decomposition of $f$ of the form (\ref{1Aaaa}).

When $0<p\leq 1,$ the Hardy martingale space $H_{p}$ is proper subspace of
Lebesque space $L_{p}.$ It is well known that for $1<p<\infty $ the space $%
H_{p}$ is nothing but $L_{p}.$

The classical inequality of Hardy type is well known in the trigonometric as
well as in the Vilenkin-Fourier analysis. Namely,
\begin{equation*}
\overset{\infty }{\underset{k=1}{\sum }}\frac{\left| \widehat{f}\left(
k\right) \right| }{k}\leq c\left\| f\right\| _{H_{1}},
\end{equation*}
where the function $f$ belongs to the Hardy space $H_{1}$ and $c$ is an
absolute constant. This was proved in the trigonometric case by Hardy and
Littlewood \cite{hl} (see also Coifman and Weiss \cite{cw}) and for Walsh
system it can be found in \cite{sws}.

Weisz \cite{We1,We4} generalized this result for Vilenkin system and proved:

\textbf{Theorem A.} Let $0<p\leq 2.$ Then there is an absolute constant $%
c_{p},$ defend only $p,$ such that
\begin{equation}
\overset{\infty }{\underset{k=1}{\sum }}\frac{\left| \widehat{f}\left(
k\right) \right| ^{p}}{k^{2-p}}\leq c_{p}\left\| f\right\| _{H_{p}},
\label{1bb}
\end{equation}
for all $f\in H_{p}.$

Paley \cite{p} proved that the Walsh- Fourier coefficients of a function $%
f\in L_{p}$ $\left( 1<p<2\right) $ satisfy the condition
\begin{equation*}
\overset{\infty }{\underset{k=1}{\sum }}\left\vert \widehat{f}\left(
2^{k}\right) \right\vert ^{2}<\infty .
\end{equation*}%
This results fails to hold $for$ $p=1.$ However, it can be verified for
functions $f\in L_{1},$ such that $f^{\ast }$ belongs $L_{1}.$ i.e $f\in
H_{1}$ (see e.g Coifman and Weiss \cite{cw})$.$

For the Vilenkin system the following theorem (see Weisz \cite{we5}) is
proved:

\textbf{Theorem B.} Let $0<p\leq 1.$ Then there is an absolute constant $%
c_{p},$ defend only $p,$ such that
\begin{equation}
\left( \underset{k=1}{\overset{\infty }{\sum }}M_{k}^{2-2/p}\underset{j=1}{%
\overset{m_{k}-1}{\sum }}\left| \widehat{f}\left( jM_{k}\right) \right|
^{2}\right) ^{1/2}\leq c_{p}\left\| f\right\| _{H_{p}},  \label{1aa}
\end{equation}
for all $f\in H_{p}.$

It is well-known that Vilenkin system forms not basis in the space $L_{1}.$
Moreover, there is a function in the dyadic Hardy space $H_{1},$ such that
the partial sums of $f$ are not bounded in $L_{1}$-norm. However, in Simon
\cite{Si3} the following strong convergence result was obtained for all $%
f\in H_{1}:$%
\begin{equation*}
\underset{n\rightarrow \infty }{\lim }\frac{1}{\log n}\overset{n}{\underset{%
k=1}{\sum }}\frac{\left\Vert S_{k}f-f\right\Vert _{1}}{k}=0,
\end{equation*}%
where $S_{k}f$ denotes the $k$-th partial sum of the Walsh-Fourier series of
$f.$ (For the trigonometric analogue see Smith \cite{sm}, for the Vilenkin
system by Gát \cite{gat1}). For the Vilenkin system Simon \cite{si1} proved:

\textbf{Theorem C.} Let $0<p<1.$ Then there is an absolute constant $c_{p},$
depends only $p,$ such that
\begin{equation}
\overset{\infty }{\underset{k=1}{\sum }}\frac{\left\Vert S_{k}f\right\Vert
_{p}^{p}}{k^{2-p}}\leq c_{p}\left\Vert f\right\Vert _{H_{p}}^{p},
\label{1cc}
\end{equation}%
for all $f\in H_{p}.$

Strong convergence theorems of two-dimensional partial sums was investigate
by Weisz \cite{We}, Goginava \cite{gg}, Gogoladze \cite{Go}, Tephnadze \cite%
{tep1}.

The main aim of this paper is to prove that the following is true:

\begin{theorem}
Let $\left\{ \Phi _{n}\right\} _{n=1}^{\infty }$ is any nondecreasing
sequence, satisfying the condition $\underset{n\rightarrow \infty }{\lim }%
\Phi _{n}=+\infty .$ Then there exists a martingale $f\in H_{p},$ such that
\begin{equation}
\text{ }\underset{k=1}{\overset{\infty }{\sum }}\frac{\left\vert \widehat{f}%
\left( k\right) \right\vert ^{p}\Phi _{k}}{k^{2-p}}=\infty ,\text{ \qquad
for }0<p\leq 2,  \label{2a}
\end{equation}%
\begin{equation}
\underset{k=1}{\overset{\infty }{\sum }}\frac{\Phi _{M_{k}}}{M_{k}^{2/p-2}}%
\underset{j=1}{\overset{m_{k}-1}{\sum }}\left\vert \widehat{f}\left(
jM_{k}\right) \right\vert ^{2}=\infty ,\text{ \qquad for }0<p\leq 1
\label{2b}
\end{equation}%
and
\end{theorem}

\begin{equation}
\text{ }\underset{k=1}{\overset{\infty }{\sum }}\frac{\left\Vert
S_{k}f\right\Vert _{L_{p,\infty }}^{p}\Phi _{k}}{k^{2-p}}=\infty ,\text{
\qquad for }0<p<1.  \label{2c}
\end{equation}

\textbf{Proof of Theorem 1.} Let $0<p\leq 2$ and $\left\{ \Phi _{n}\right\}
_{n=1}^{\infty }$ is any nondecreasing, nonnegative sequence, satisfying
condition
\begin{equation*}
\underset{n\rightarrow \infty }{\lim }\Phi _{n}=\infty .
\end{equation*}%
For this function $\Phi \left( n\right) ,$ there exists an increasing
sequence $\left\{ \alpha _{k}\geq 2:k\in \mathbb{N}_{+}\right\} $ of the
positive integers such that:

\begin{equation}
\sum_{k=1}^{\infty }\frac{1}{\Phi _{M_{_{\alpha _{k}}}}^{p/4}}<\infty .
\label{2}
\end{equation}

Let \qquad
\begin{equation*}
f^{\left( A\right) }\left( x\right) :=\sum_{\left\{ k;\text{ }\alpha
_{k}<A\right\} }\lambda _{k}a_{k}\left( x\right) ,
\end{equation*}%
where
\begin{equation*}
\lambda _{k}=\frac{1}{\Phi _{M_{_{\alpha _{k}}}}^{1/4}}
\end{equation*}%
and

\begin{equation*}
a_{k}\left( x\right) =\frac{M_{\alpha _{k}}^{1/p-1}}{M}\left( D_{M_{\alpha
_{k}+1}}\left( x\right) -D_{M_{_{\alpha _{k}}}}\left( x\right) \right) ,
\end{equation*}%
where $M=\sup_{n\in \mathbb{N}}m_{n}.$

It is easy to show that the martingale $\,f=\left( f^{\left( 1\right)
},f^{\left( 2\right) },...,f^{\left( A\right) },...\right) \in H_{p}.$

Indeed, since

\begin{equation}
S_{M_{A}}\left( a_{k}\left( x\right) \right) =\left\{
\begin{array}{l}
a_{k}\left( x\right) \text{ \qquad }\alpha _{k}<A \\
0\text{ \qquad }\alpha _{k}\geq A,%
\end{array}%
\right.  \label{4}
\end{equation}

\begin{eqnarray*}
\text{supp}(a_{k}) &=&I_{\alpha _{k}}, \\
\int_{I_{\alpha _{k}}}a_{k}d\mu &=&0,
\end{eqnarray*}

and

\begin{eqnarray*}
\left\Vert a_{k}\right\Vert _{\infty } &\leq &\frac{M_{\alpha _{k}}^{1/p-1}}{%
M}M_{\alpha _{k}+1} \\
&\leq &M_{\alpha _{k}}^{1/p}=\mu (\text{supp }a_{k})^{-1/p},
\end{eqnarray*}%
if we apply Theorem W and (\ref{2}) we conclude that $f\in H_{p}.$

It is easy to show that

\begin{equation}
\widehat{f}(j)=\left\{
\begin{array}{l}
\frac{1}{M}\frac{M_{\alpha _{k}}^{1/p-1}}{\Phi _{M_{_{\alpha _{k}}}}^{1/4}}%
,\,\qquad \,\text{ if \thinspace \thinspace }j\in \left\{ M_{\alpha
_{k}},...,\text{ ~}M_{\alpha _{k}+1}-1\right\} ,\text{ }k=1,2... \\
0,\text{ }\qquad \text{\thinspace \thinspace \thinspace if \thinspace
\thinspace \thinspace }j\notin \bigcup\limits_{k=1}^{\infty }\left\{
M_{\alpha _{k}},...,\text{ ~}M_{\alpha _{k}+1}-1\right\} .\text{ }%
\end{array}%
\right.  \label{5}
\end{equation}%
First we prove equality (\ref{2a}). Using (\ref{5}) we can
\begin{eqnarray*}
&&\underset{l=1}{\overset{M_{\alpha _{k}+1}-1}{\sum }}\frac{\left\vert
\widehat{f}\left( l\right) \right\vert ^{p}\Phi _{l}}{l^{2-p}} \\
&=&\underset{n=1}{\overset{k}{\sum }}\underset{l=M_{\alpha _{n}}}{\overset{%
M_{\alpha _{n}+1}-1}{\sum }}\frac{\left\vert \widehat{f}\left( l\right)
\right\vert ^{p}\Phi _{l}}{l^{2-p}} \\
&\geq &\underset{l=M_{\alpha _{k}}}{\overset{M_{\alpha _{k}+1}-1}{\sum }}%
\frac{\left\vert \widehat{f}\left( l\right) \right\vert ^{p}\Phi _{l}}{%
l^{2-p}} \\
&\geq &c\Phi _{M_{\alpha _{k}}}\underset{l=M_{\alpha _{k}}}{\overset{%
M_{\alpha _{k}+1}-1}{\sum }}\frac{\left\vert \widehat{f}\left( l\right)
\right\vert ^{p}}{l^{2-p}} \\
&\geq &c\Phi _{M_{\alpha _{k}}}\frac{M_{\alpha _{k}}^{1-p}}{\Phi _{M_{\alpha
_{k}}}^{p/4}}\underset{l=M_{\alpha _{k}}}{\overset{M_{\alpha _{k}+1}-1}{\sum
}}\frac{1}{l^{2-p}} \\
&\geq &c\Phi _{M_{\alpha _{k}}}^{1/2}M_{\alpha _{k}}^{1-p}\underset{%
l=M_{\alpha _{k}}}{\overset{M_{\alpha _{k}+1}-1}{\sum }}\frac{1}{M_{\alpha
_{k}+1}^{2-p}} \\
&\geq &c\Phi _{M_{\alpha _{k}}}^{1/2}M_{\alpha _{k}}^{1-p}\frac{1}{M_{\alpha
_{k}+1}^{1-p}} \\
&\geq &c\Phi _{M_{\alpha _{k}}}^{1/2}\rightarrow \infty ,\text{ \qquad when }%
k\rightarrow \infty .
\end{eqnarray*}

Next we prove equality (\ref{2b}).\textbf{\ }Let $0<p\leq 1.$ Using (\ref{5}%
) we get%
\begin{eqnarray*}
&&\underset{l=1}{\overset{k}{\sum }}M_{\alpha _{l}}^{2-2/p}\Phi _{M_{\alpha
_{l}}}\underset{j=1}{\overset{m_{\alpha _{l}}-1}{\sum }}\left\vert \widehat{f%
}\left( jM_{\alpha _{l}}\right) \right\vert ^{2} \\
&\geq &M_{\alpha _{k}}^{2-2/p}\Phi _{M_{\alpha _{k}}}\underset{j=1}{\overset{%
m_{\alpha _{k}}-1}{\sum }}\left\vert \widehat{f}\left( jM_{\alpha
_{k}}\right) \right\vert ^{2} \\
&\geq &cM_{\alpha _{k}}^{2-2/p}\Phi _{M_{\alpha _{k}}}\underset{j=1}{\overset%
{m_{\alpha _{k}}-1}{\sum }}\frac{M_{\alpha _{k}}^{2/p-2}}{\Phi _{M_{\alpha
_{k}}}^{1/2}}
\end{eqnarray*}%
\begin{equation*}
\geq c\Phi _{M_{\alpha _{k}}}^{1/2}\rightarrow \infty ,\text{ \qquad when }%
k\rightarrow \infty .
\end{equation*}

Finally we prove equality (\ref{2c}).\textbf{\ }Let\textbf{\ }$0<p<1$ and $%
M_{\alpha _{k}}\leq j<M_{\alpha _{k}+1}$. From (\ref{5}) we have

\begin{eqnarray*}
&&S_{j}f\left( x\right) =\sum_{l=0}^{M_{\alpha _{k-1}+1}-1}\widehat{f}%
(l)\psi _{l}\left( x\right) \\
&&+\sum_{l=M_{\alpha _{k}}}^{j-1}\widehat{f}(l)\psi _{l}\left( x\right) \\
&=&\sum_{\eta =0}^{k-1}\sum_{v=M_{\alpha _{\eta }}}^{M_{\alpha _{\eta }+1}-1}%
\widehat{f}(v)\psi _{v}\left( x\right) \\
&&+\sum_{v=M_{\alpha _{k}}}^{j-1}\widehat{f}(v)\psi _{v}\left( x\right) \\
&=&\sum_{\eta =0}^{k-1}\sum_{v=M_{\alpha _{\eta }}}^{M_{\alpha _{\eta }+1}-1}%
\frac{1}{M}\frac{M_{\alpha _{\eta }}^{1/p-1}}{\Phi _{M_{\alpha _{\eta
}}}^{1/4}}\psi _{v}\left( x\right) \\
&&+\sum_{v=M_{\alpha _{k}}}^{j-1}\frac{1}{M}\frac{M_{\alpha _{k}}^{1/p-1}}{%
\Phi _{M_{\alpha _{k}}}^{1/4}}\psi _{v}\left( x\right) \\
&=&\sum_{\eta =0}^{k-1}\frac{1}{M}\frac{M_{\alpha _{\eta }}^{1/p-1}}{\Phi
_{M_{\alpha _{\eta }}}^{1/4}}\left( D_{M_{_{\alpha _{\eta }+1}}}\left(
x\right) -D_{M_{_{\alpha _{\eta }}}}\left( x\right) \right) \\
&&+\frac{1}{M}\frac{M_{\alpha _{k}}^{1/p-1}}{\Phi _{M_{\alpha _{k}}}^{1/4}}%
\left( D_{_{j}}\left( x\right) -D_{M_{_{\alpha _{k}}}}\left( x\right) \right)
\\
&=&I+II.
\end{eqnarray*}

Let $j\in \mathbb{N}_{n_{0}}$ and $x\in G_{m}\backslash I_{1}.$ Since $%
j-M_{\alpha _{k}}\in \mathbb{N}_{n_{0}}$ and
\begin{equation*}
D_{j+M_{\alpha _{k}}}\left( x\right) =D_{M_{\alpha _{k}}}\left( x\right)
+\psi _{M_{\alpha _{k}}}\left( x\right) D_{j}\left( x\right) ,\text{ when }%
\,\,j<M_{\alpha _{k}},
\end{equation*}%
combining (\ref{3}) and (\ref{3aa}) we can write

\begin{eqnarray}
\left\vert II\right\vert &=&\frac{1}{M}\frac{M_{\alpha _{k}}^{1/p-1}}{\Phi
_{M_{\alpha _{k}}}^{1/4}}\left\vert \psi _{M_{\alpha
_{k}}}D_{_{j-M_{_{\alpha _{k}}}}}\left( x\right) \right\vert  \label{13a} \\
&=&\frac{1}{M}\frac{M_{\alpha _{k}}^{1/p-1}}{\Phi _{M_{\alpha _{k}}}^{1/4}}%
\left\vert \psi _{M_{\alpha _{k}}}\left( x\right) \psi _{j-M_{_{\alpha
_{k}}}}\left( x\right) r_{0}^{m_{0}-1}\left( x\right) D_{_{1}}\left(
x\right) \right\vert  \notag \\
&=&\frac{1}{M}\frac{M_{\alpha _{k}}^{1/p-1}}{\Phi _{M_{\alpha _{k}}}^{1/4}}.
\notag
\end{eqnarray}%
Applying (\ref{3}) and condition $\alpha _{n}\geq 2$ $\left( n\in \mathbb{N}%
\right) $ for $I$ we have

\begin{equation}
I=0,\text{ \qquad for }x\in G_{m}\backslash I_{1}.  \label{13b}
\end{equation}%
It follows that
\begin{equation*}
\left\vert S_{j}f\left( x\right) \right\vert =\left\vert II\right\vert =%
\frac{1}{M}\frac{M_{\alpha _{k}}^{1/p-1}}{\Phi _{M_{\alpha _{k}}}^{1/4}},%
\text{ \qquad for }x\in G_{m}\backslash I_{1}.
\end{equation*}%
Hence

\begin{eqnarray}
&&\left\Vert S_{j}\left( f\left( x\right) \right) \right\Vert _{L_{p,\infty
}}  \label{13} \\
&\geq &\frac{1}{2M}\frac{M_{\alpha _{k}}^{1/p-1}}{\Phi _{M_{\alpha
_{k}}}^{1/4}}\mu \left( x\in G_{m}:\left\vert S_{j}\left( f\left( x\right)
\right) \right\vert >\frac{1}{2M}\frac{M_{\alpha _{k}}^{1/p-1}}{\Phi
_{M_{\alpha _{k}}}^{1/4}}\right) ^{1/p}  \notag \\
&\geq &\frac{1}{2M}\frac{M_{\alpha _{k}}^{1/p-1}}{\Phi _{M_{\alpha
_{k}}}^{1/4}}\mu \left( x\in G_{m}\backslash I_{1}:\left\vert S_{j}\left(
f\left( x\right) \right) \right\vert >\frac{1}{2M}\frac{M_{\alpha
_{k}}^{1/p-1}}{\Phi _{M_{\alpha _{k}}}^{1/4}}\right) ^{1/p}  \notag \\
&=&\frac{1}{2M}\frac{M_{\alpha _{k}}^{1/p-1}}{\Phi _{M_{\alpha _{k}}}^{1/4}}%
\left\vert G_{m}\backslash I_{1}\right\vert  \notag \\
&\geq &\frac{cM_{\alpha _{k}}^{1/p-1}}{\Phi _{M_{\alpha _{k}}}^{1/4}}.
\notag
\end{eqnarray}

Combining (\ref{1a}) and (\ref{13}) we have

\begin{eqnarray*}
&&\underset{j=1}{\overset{M_{\alpha _{k}+1}-1}{\sum }}\frac{\left\Vert
S_{j}\left( f\left( x\right) \right) \right\Vert _{L_{p,\infty }}^{p}\Phi
_{j}}{j^{2-p}} \\
&\geq &\underset{j=M_{\alpha _{k}}}{\overset{M_{\alpha _{k}+1}-1}{\sum }}%
\frac{\left\Vert S_{j}\left( f\left( x\right) \right) \right\Vert
_{L_{p,\infty }}^{p}\Phi _{j}}{j^{2-p}} \\
&\geq &\Phi _{M_{\alpha _{k}}}\underset{\left\{ n:M_{k}\leq n\leq M_{k+1},%
\text{ }n\in \mathbb{N}_{n_{0}}\right\} }{\sum }\frac{\left\Vert S_{j}\left(
f\left( x\right) \right) \right\Vert _{L_{p,\infty }}^{p}}{j^{2-p}} \\
&\geq &c\Phi _{M_{\alpha _{k}}}\frac{M_{\alpha _{k}}^{1-p}}{\Phi _{M_{\alpha
_{k}}}^{p/4}}\underset{\left\{ n:M_{k}\leq n\leq M_{k+1},\text{ }n\in
\mathbb{N}_{n_{0}}\right\} }{\sum }\frac{1}{j^{2-p}} \\
&\geq &c\Phi _{M_{\alpha _{k}}}^{3/4}M_{\alpha _{k}}^{1-p}\underset{\left\{
n:M_{k}\leq n\leq M_{k+1},\text{ }n\in \mathbb{N}_{n_{0}}\right\} }{\sum }%
\frac{1}{M_{\alpha _{k}+1}^{2-p}}
\end{eqnarray*}%
\begin{eqnarray*}
&\geq &c\frac{\Phi _{M_{\alpha _{k}}}^{3/4}}{M_{\alpha _{k}+1}}\underset{%
\left\{ n:M_{k}\leq n\leq M_{k+1},\text{ }n\in \mathbb{N}_{n_{0}}\right\} }{%
\sum }1 \\
&\geq &c\Phi _{M_{\alpha _{k}}}^{3/4}\rightarrow \infty ,\text{ \qquad when }%
k\rightarrow \infty .
\end{eqnarray*}

Theorem 1 is proved.

\end{document}